\documentclass[11pt,fleqn]{article}

\usepackage{latexsym}

\usepackage{amsmath}
\usepackage{amsthm}
\usepackage{amssymb}
\usepackage{amsfonts}

\newcommand{\be}{\begin{equation} } 
\newcommand{\ee}{\end{equation} \par \noindent}

\newcommand{\rf}[1]{(\ref{#1})}

\newcommand{\scal}[2]{\mbox{$\langle #1 \! \mid #2 \rangle $}} 

\newcommand{\ba}{\begin{array}}
\newcommand{\ea}{\end{array}}
\newcommand{\const}{{\rm const}}

\newcommand{\e}{{\mathbf e}}

\newcommand{\R}{{\mathbb R}}
\newcommand{\Z}{{\mathbb Z}}

\newcommand{\T}{{\mathbb T}}

\newcommand{\dis}{\displaystyle}

\newcommand{\rr}{{\vec r}}
\newcommand{\nn}{{\vec n}}
\newcommand{\ep}{\varepsilon}

\newcommand{\ema}{\ea \right)}
\newcommand{\mm}{\left( \ba{rr} }

\newtheorem{prop}{Proposition}
\newtheorem{Th}{Theorem}

\newtheorem{cor}{Corollary}
\newtheorem{Def}{Definition}

\newenvironment{Proof}{\par \vspace{2ex} \par
\noindent \small {\it Proof:}}{\hfill $\Box$ 
\vspace{2ex} \par }

\begin{document}

\title{\bf 
Pseudospherical surfaces on time scales: \\ a geometric definition and the spectral approach}

\author{{\bf
Jan L.\ Cie\'sli\'nski\thanks{
E-mail: \tt janek\,@\,alpha.uwb.edu.pl}}
\\ {\footnotesize Uniwersytet w Bia\l ymstoku, Instytut Fizyki Teoretycznej}
\\ {\footnotesize ul.\ Lipowa 41, 15-424  Bia\l ystok, Poland} 
}

\maketitle

\begin{abstract} 
We define and discuss the notion of pseudospherical surfaces in asymptotic coordinates on time scales. Thus we extend well known notions of discrete pseudospherical surfaces and smooth pseudosperical surfaces on more exotic domains (e.g, the Cantor set). In particular, we present a new expression for the discrete Gaussian curvature which turns out to be valid for asymptotic nets on any time scale. We show that asymptotic Chebyshev nets on an arbitrary time scale have constant negative Gaussian curvature. We present also the quaternion-valued spectral problem (the Lax pair) and the Darboux-B\"acklund transformation for pseudospherical surfaces (in asymptotic coordinates) on arbitrary time scales. 
\end{abstract}

\par \vspace{0.5cm} \par
\par \vspace{0.5cm} \par
\par \vspace{0.5cm} \par

\section{Introduction}

A time scale (or a measure chain) is an arbitrary nonempty closed subset of the real numbers \cite{Hi}. Typical examples are $\R$,  $\Z$, any unions of isolated points and closed intervals, and, finally, discrete sets containing all acumulation points (like the Cantor set). The time scales were introduced in order to unify differential and difference calculus \cite{Hi,Hi2}. 
Partial differentiation, tangent lines and tangent planes on time scales have been introduced recently \cite{BG}. In this paper we suggest how to extend the differentiation also on Lie groups. The case of the $SU(2)$ group is discussed in detail.

The difference geometry \cite{Sau} is a discrete analogue of the differential geometry. In the last years one can observe a fast development of the integrable difference geometry (see, for instance, \cite{BP-rev,BS,Ci-fam,CDS,Do,DS-nets,WKS}) closely related to the classical differential geometry \cite{DSM,Eis1}. It is interesting that in the discrete case one recovers explicit constructions and transformations known in the continuous case (e.g., Darboux, B\"acklund, Ribaucour, Laplace and Jonas transformations, soliton and finite-gap solutions etc.). 
A natural idea is to unify the difference and differential geometries and to formulate the integrable geometry on time scales.

In this paper we propose such formulation for pseudospherical immersions (surfaces of constant negative Gaussian curvature). The discrete pseudospherical surfaces have been introduced a long time ago \cite{Sau2,Wun}, and studied intensively in the last years \cite{BP-pseudo}. 
The idea to extend the notion of pseudospherical surfaces on arbitrary time scales first appeared in \cite{Sw-mgr}. However, throughout that work 
there was assumed  that all points are isolated (the discrete case). The discrete Gaussian curvature and the B\"acklund transformation were not considered at all. In the present paper we formulate a natural geometric definition of pseudospherical surfaces (more precisely: asymptotic Chebyshev nets) on time scales and present the associated spectral problem (the Lax pair) and the Darboux-B\"acklund transformation. Thus the discrete, continuous and other cases are first described in a unified framework.

\section{Differentiation on time scales}

This section collects basic notions and results concerning the differential calculus on time scales, compare \cite{BG}. To avoid some unimportant complications we confine ourselves to time scales which are not bounded neither from above nor from below.

\begin{Def}[\cite{Hi}] Let a time scale $\T$ is given. The maps $\sigma : \T \rightarrow \T$ and $\rho : \T \rightarrow \T$, defined by 
\be  
\sigma (u) := \inf \{ v \in \T : v > u \}  \ ,  \qquad
\rho (u) := \sup \{ v \in \T : v < u \}  \ ,
\ee
are called jump operator and backward jump operator, respectively.
\end{Def}

\begin{Def}[\cite{Hi}] \label{class} A point $u \in \T$ is said to be right-scattered (if $\sigma (u) > u$) or right-dense (if $\sigma (u) = u$), 
 left-scattered (if $\rho (u) < u$) or left-dense (if $\rho (u) = u$), and isolated if \ $\rho (u) < u < \sigma (u)$.
\end{Def}

\begin{Def}[\cite{BG}]
The delta derivative of a continuous function $f$ is defined as
\be
\frac{\partial f (t) }{\Delta t} = \lim_{\stackrel{s \rightarrow t}{s \neq \sigma (t)}} \frac{ f ( \sigma (t) ) - f ( s )}{\sigma (t) - s} \ ,
\ee
and the nabla derivative is defined by
\be
\frac{\partial f (t) }{\nabla t} = \lim_{\stackrel{s \rightarrow t}{s \neq \rho (t)}} \frac{ f ( \rho (t) ) - f ( s )}{\rho (t) - s} \ .
\ee
\end{Def}

In this paper we focus on functions defined on two-dimensional time scales, i.e., on $\T_1 \times \T_2 $, where $\T_1, \T_2$ are given time scales. 
The extension on $n$-dimensional time scales is usually straightforward. We denote: 
\be  \ba{l}
t \equiv (t_1, t_2) \in \T_1 \times \T_2 \ ,  \\[2ex]
\sigma_1 (t) = ( \sigma (t_1) , t_2)  \ , \quad \sigma_2 (t) = (t_1, \sigma (t_2)) \ , \\[2ex]
\rho_1 (t) = ( \rho (t_1) , t_2)  \ , \quad \rho_2 (t) = (t_1, \rho (t_2)) \ .
\ea \ee

Obviously, the notions from Definition~\ref{class} can be defined independently for each variable. For example, a point can be right-scattered in the first variable, and right-dense in the second variable, 
shortly: 1-right-scattered and 2-right-dense. 

We stress that throughout this paper $\sigma_1$ and $\sigma_2$ usually denote jump operators, unless stated otherwise (only in few places in the text we mention Pauli sigma matrices denoted by ${\pmb \sigma}_1$, ${\pmb \sigma}_2$, ${\pmb \sigma}_3$).

In the discrete case (e.g., $\T_1 = \T_2 = \Z$)  we have $\sigma_j (u) = T_j u$ and $\rho_j (u) = T_j^{-1} u$, where $T_1$, $T_2$ are usual shift operators. Therefore delta and nabla differentiation can be associated with forward and backward data, respectively \cite{Do}.

\begin{Def}[\cite{BG}]
The partial delta derivative of a continuous function $f$ is defined as
\be
\frac{\partial f (t) }{\Delta t_j} = \lim_{\stackrel{s_j \rightarrow t_j}{s_j \neq \sigma (t_j)}} \frac{ f ( \sigma_j (t) ) - f ( s )}{\sigma (t_j) - s_j} \ .
\ee
The definition of the partial nabla derivative is analogical.
\end{Def}

\begin{prop}[\cite{BG}]
If the mixed partial delta derivatives exist in a neighbourhood of $t_0 \in \T_1 \times \T_2$  and are continuous at $t = t_0$, then 
\[
   \frac{\partial^2 f (t_0) }{\Delta t_1  \Delta t_2 } = \frac{\partial^2 f (t_0) }{\Delta t_2 \Delta t_1 } \ .
\]
\end{prop}

In the continuous case (e.g., $\T_1 = \T_2 = \R$) the delta derivative coincides with the right-hand derivative, while the nabla derivative coincides with the left-hand derivative. Note that all results and definitions in terms of delta derivatives have their nabla derivatives analogues. 

In the continuous case the differentiability implies the existence of the tangent plane. The delta differentiability does not have this important property. We need a stronger notion: the complete delta differentiability.  

\begin{Def}[\cite{BG}]  \label{complete}
 We say that a function $f : \T \rightarrow \R$ is completely delta differentiable at a point $t_0 \in \T$, if there exist a number $A$ such that 
\[  \ba{l}
f (t) - f(t_0) = A (t - t_0) + (t - t_0) \ \alpha (t_0, t) \ , \\[2ex]
f(t) - f(\sigma (t_0)) = A ( t - \sigma (t_0) ) + (t - \sigma (t_0)) \ \beta (t_0, t)
\ ,
\ea  \]
where \ $\alpha (t_0, t_0) = 0$, $\beta (t_0, t_0) = 0$, $\dis \lim_{t\rightarrow t_0} \alpha (t_0, t) = 0$, and \ $\dis  \lim_{t\rightarrow t_0} \beta (t_0, t) = 0$. 
\end{Def}

\begin{prop}[\cite{BG}]  \label{tanline}
If the function $f$ is completely delta differentiable at $t_0$, then the graph of this function has the uniquely determined delta tangent line at the point $P_0 = (t_0, f(t_0))$ specified by the equation
\[
y - f(t_0) = \frac{\partial f (t_0) }{\Delta t} (x - t_0) 
\] 
\end{prop}

If $P_0$ is an isolated point of the curve $\Gamma$ (hence $P_0 \neq P_0^\sigma$), then 
the delta tangent line to $\Gamma$ at $P_0$ coincides with the unique line through the points $P_0$ and $P_0^\sigma$.

The definition of the complete delta differentiability in two-dimensional case is similar to Definition~\ref{complete} (for details, see \cite{BG}, Definition 2.1). Instead of this definition we present here an important sufficient condition for delta differentiability.

\begin{prop}[\cite{BG}]
Let $f : \T_1 \times \T_2 \rightarrow \R$ is continuous and has first order partial delta derivatives in a neighbourhood of $t_0$. If these derivatives are continuous at $t_0$, then $f$ is completely delta differentiable at $t_0$.
\end{prop}

If $P_0 \neq P_0^{\sigma_1}$ and $P_0^{\sigma_2} \neq P_0$ (hence also $P_0^{\sigma_1} \neq P_0^{\sigma_2}$), then the delta tangent plane to the surface $S$ at $P_0$ (if exists) coincides with the unique plane through $P_0, P_0^{\sigma_1}$ and $P_0^{\sigma_2}$.

\begin{prop}[\cite{BG}]   \label{tanplane}
If the function $f : \T_1 \times \T_2 \rightarrow {\mathbb R}$ is completely delta differentiable at \ $t_0 = (t_{01}, t_{02})$, then the surface represented by this function has the uniquely determined delta tangent plane at the point $P_0 = (t_{01}, t_{02}, f(t_0))$ specified by the equation 
\be
z = f(t_0) + \frac{\partial f (t_0)}{\Delta t_1} (x - t_{01}) + \frac{\partial f (t_0)}{\Delta t_2} (y - t_{02}) 
\ee
where $(x,y,z)$ is the current point of the plane.
\end{prop}

In the following sections of this paper we define pseudospherical surfaces on time scales in terms of delta derivatives. In order to simplify the notation the delta derivatives will be denoted by 
\be  \label{D}
D_j f \equiv \frac{\partial f (t) }{\Delta t_j} \ .
\ee
Propositions~\ref{tanline} and \ref{tanplane} show that in  geometrical contexts the complete delta differentiability, which guarantees the existence of tangent lines and tangent planes, is more useful than the delta differentiability.

\section{Differentiation of $SU(2)$-valued functions on time scales}

Analytic approaches to pseudospherical surfaces usually involve the Lie group $SU(2)$, Lie algebra $su(2)$ and quaternions \cite{BP-rev,MeS,PM,Sym}. Therefore it is important to extend the notion of the delta derivative on Lie groups.  

Given a function $f : {\mathbb T} \rightarrow M$, where $M$ is a submanifold, we can define the delta derivative of $f$ in a quite natural way. If $t$ is right-dense, then we compute the tangent vector in the point $t$ just repeating the standard procedure, well known in the case ${\mathbb T} = {\mathbb R}$. 
If $t$ is right-scattered, then we join $f(t)$ and $f(\sigma(t))$ by the shortest geodesic. The delta derivative is defined as the vector tangent to this geodesic. 
If $M = G$ is a Lie group, then we may map the tangent vector into the coresponding Lie algebra $g$. 
The length of this vector is $\delta/\ep$, where $\ep = \sigma(t) - t$ and $\delta$ is the lenght of the geodesic between $t$ and $\sigma(t)$. 
 
If $M$ is immersed in an ambient Euclidean space, then one can define the delta derivative in another way, considering geodesics (straight lines) in the ambient space instead of geodesics on $M$. Both definitions yield the same results for right-dense points, but for right-scattered points we get two different definitions of the delta derivative (even after projection onto the corrsponding tangent space).  
In the general case these ideas will be developed elsewhere. Here we confine ourselves to the Lie group $SU(2)$.

The Lie group $SU(2)$ is defined as $\{ \Phi : \Phi^{-1} = \Phi^\dagger, \ \det \Phi = 1 \}$. Any element $\Phi \in SU(2)$ can be parameterized as
\be   \label{Su2}
  \Phi = \mm a & b \\ - \bar b & \bar a \ema  \ , \quad |a|^2 + |b|^2 = 1 \ .
\ee
Therefore 
\be
\Phi = {\rm Re} a  -  \e_1 {\rm Re} b  -  \e_2 {\rm Im\, } b  -  \e_3 {\rm Im\, } a \ , 
\ee
where $\e_j = - i {\pmb \sigma}_j$ ($j=1,2,3$) and ${\pmb \sigma}_j$ are standard Pauli matrices. The following properties are satisfied: 
\be \label{ijk}
\e_1^2 = \e_2^2 = \e_3^2 = - 1 \ , \qquad \e_j \e_k = - \e_k \e_j \quad (j\neq k) \ ,
\ee
\be  \label{123}
\e_1 \e_2 = \e_3 \ , \quad \e_2 \e_3 = \e_1 \ , \quad 
\e_3 \e_1 = \e_2 \ , 
\ee
\be  \label{edag}
\e_j^\dagger = - \e_j \ \ (j = 1, 2, 3) \ .
\ee
Therefore the space spanned by $1, \e_1 , \e_2 , \e_3$ can be identified with quaternions ${\mathbb H}$. The standard Euclidean structure is defined by the following scalar product 
\be
   \scal{A}{B}  = \frac{1}{2} {\rm Tr} (A B^\dagger)  \ , \qquad A, B \in {\mathbb H} \ .
\ee
Then the basis $1, \e_1 , \e_2 , \e_3$ is orthonormal. The space of imaginary (or pure) quaternions, ${\rm Im\, } {\mathbb H}$, is spanned by $\e_1, \e_2, \e_3$.

The condition $|a|^2 + |b|^2 = 1$ means exactly that 
$\Phi$ given by \rf{Su2} is a unit vector. 
Hence we have the well known conclusion that the Lie group $SU(2)$ can be identified with the sphere $S^3 \subset \mathbb H$. The Lie algebra $su(2)$ coincides with pure quaternions ${\rm Im\, } {\mathbb H}$.

Following the general outline given above we are going to define two delta derivatives, denoted by ${\cal D}_j$ and $D_j$, respectively. In the continuous case ($j$-right-dense points) ${\cal D}_j = D_j$ and
\be
U_j := (D_j \Phi ) \Phi^{-1} 
\ee
takes values in the Lie algebra $su(2)$. 
In the discrete case ($j$-right-scattered points) the situation is more complicated.

Geometrically, the derivative ${\cal D}_j \Phi$ in the discrete case is tangent to the sphere $S^3$ at $\Phi$ and $|{\cal D}_j \Phi|$  is the length of the coresponding arc. Therefore, after elementary geometric considerations, 
\be
{\cal D}_j \Phi = \frac{ (T_j (\Phi) - \Phi \cos\delta ) \delta }{\ep \sin\delta} \ , \qquad \cos\delta := \scal{T_j \Phi}{\Phi} \ . 
\ee
Note that 
\be
T_j \Phi = \exp ( {\mathbf  u}_j \delta ) \Phi \ , \qquad {\mathbf  u}_j := \frac{\ep}{\delta} ( {\cal D}_j \Phi ) \Phi^{-1} \ ,
\ee
and ${\mathbf  u}_j$ is a unit vector from $su(2)$.
The derivative $D_j \Phi$ can be identified  with the secant joining $\Phi$ and $T_j \Phi$ (in the space ${\mathbb H}$):
\be  \label{DFi}
D_j \Phi = \frac{T_j \Phi - \Phi}{\ep} \ .
\ee
Now $(D_j \Phi) \Phi^{-1}$ is, in general, outside ${\rm Im\, } {\mathbb H}$. 
Therefore, it is convenient to define a projection $\Pi : {\mathbb H} \rightarrow { \rm Im\, } {\mathbb H}$ 
\be  \label{Pi}
   \Pi (A_0 + A_1 \e_1 + A_2 \e_2 + A_3 \e_3) := A_1 \e_1 + A_2 \e_2 + A_3 \e_3 \ ,
\ee
projecting a quaternion $A$ into its imaginary (or traceless) part. 
One can check that
\be
\Pi ( ( D_j \Phi) \Phi^{-1} ) = \frac{\delta}{\sin\delta}  ( {\cal D}_j \Phi ) \Phi^{-1}   \ .
\ee

Throughout this paper we will use only the derivative $D_j$, defined by \rf{DFi}, but applied not only to elements of $SU(2)$ but to any $\Psi \in {\mathbb H}$. Note that 
\be
 \Psi = \mm a & b \\ - \bar b & \bar a \ema  \quad \Longrightarrow \quad 
 D_j \Psi = \mm D_j a & D_j b \\  - D_j \bar b & D_j \bar a \ema 
\ee
and the following rules of differentiation hold 
\be  \ba{l}  \label{Leib}
   D_j ( A B ) = ( D_j A ) B + \sigma_j (A) D_j B \ ,    \\[2ex]
D_j (\Psi^{-1}) = - \sigma_j (\Psi) (D_j \Psi) \Psi^{-1} \ ,
\ea \ee
where   $A, B, \Psi \in {\mathbb H}$. 
Therefore $D_j$ is  convenient in calculations and turned out to be sufficient for our purposes (see Section~\ref{Lax}, where the spectral approach for pseudoshperical immersions is presented). A formulation of another spectral approach  based on the more geometric derivative ${\cal D}_j$ is an open problem. It would be interesting to check the equivalence of both approaches.

\section{Smooth and discrete pseudospherical surfaces}
\label{pseudo-review}

Pseudospherical surfaces, i.e., surfaces (immersions) of constant negative Gaussian curvature have been studied intensively since the middle of the XIX century, starting from 1839 \cite{Min}. The famous transformations found by Bianchi, Lie and B\"acklund turned out to be milestones both in differential geometry and in the soliton theory. 
Old and recent results concerning pseudospherical surfaces, including a lot of orignal references, are collected and reported for instance in \cite{Eis1,Hey,Ov,PM}, see also \cite{Ampol}. 

Let us consider a surface immersed in $\R^3$ explicitly described by a position vector $\rr = \rr (s, t)$ (we assume that this function is sufficiently smooth). We denote the normal vector by $\nn$ and define the so called fundamental forms:
\be \ba{l}
  I := d \rr \cdot d \rr = E ds^2 + 2 F ds\, dt + G dt^2 \ , \\[2ex]
 II := - d \rr \cdot d \nn = L ds^2 + 2 M ds\, dt + N dt^2 \ ,
\ea \ee
where the center dot denotes the standard scalar product in $\R^3$ and $E, F, G$, $L, M, N$ are real functions of $s, t$. These functions have to satisfy nonlinear equations known as Gauss-Peterson-Codazzi equations. 
The Gaussian curvature can be conveniently expressed as follows
\be  \label{Kcont}
   K = \frac{\det (II)}{\det (I)} =  \frac{(\rr,_1 \cdot \nn,_1) (\rr,_2 \cdot \nn,_2) - (\rr,_1 \cdot \nn,_2 )(\rr,_2 \cdot \nn,_1)}{(\rr,_1 \cdot \rr,_1) (\rr,_2 \cdot \rr,_2) - (\rr,_1 \cdot \rr,_2)^2} \ ,
\ee
where $\rr,_1 := \partial \rr/\partial t$, $\rr,_2 := \partial \rr/\partial s$, etc.

\begin{Def}
Coordinates $s,t$ are called Chebyshev coordinates if the first fundamental form is given by $I = ds^2 + 2 \cos\phi ds\,dt + dt^2$, i.e., 
\be  \label{Cheb}
E \equiv \rr,_1 \cdot \rr,_1 = 1 \ , \quad G \equiv \rr,_2 \cdot \rr,_2 = 1  \ , \quad F \equiv \rr_1 \cdot \rr_2 =  \cos\phi \ .
\ee
If a less restrictive conditons hold: 
\be
 E,_2 = 0 \ , \quad G,_1 = 0 \ , 
\ee
then $s, t$ are called weak Chebyshev coordinates.
\end{Def}

Any weak Chebyshev coordinates $s,t$ can be transformed (at least locally)  into Chebyshev coordinates $\tilde s, \tilde t$ by an appropriate change of variables  $\tilde s = g (s)$, $\tilde t = f (t)$.

\begin{Def}
Coordinates $s,t$ are called asymptotic if the second fundamental form is given by $II = 2 M dt\,ds$, i.e., 
\be  \label{asym}
\rr,_1 \cdot \nn,_1 = \rr,_2 \cdot \nn,_2 = 0 \ , \quad \rr,_1 \cdot \nn,_2 = \rr,_2 \cdot \nn,_1 = - M \ .
\ee
\end{Def}

\begin{prop}  \label{czeb}
Asymptotic lines on a surface admit parameterization by Chebyshev coordinates 
if and only if the surface  has a constant negative Gaussian curvature. 
In this case the Gaussian curvature is given by
\be  \label{Kcontas}
   K = \frac{ - (\rr,_1 \cdot \nn,_2 )(\rr,_2 \cdot \nn,_1)}{(\rr,_1 \cdot \rr,_1) (\rr,_2 \cdot \rr,_2)  - (\rr,_1 \cdot \rr,_2)^2} = - \left( \frac{  M}{\sqrt{E} \sqrt{G} \sin\phi} \right)^2  \ ,
\ee
where $\phi$ is the angle between $\rr,_1$ and $\rr,_2$. 
\end{prop}

\vspace{2ex}

\noindent {\bf Discrete surfaces} (discrete immersions) are  defined as maps 
$$ \rr: \ep_1 \Z \times \ep_2 \Z \ni (\ep_1  m, \ep_2 n) \rightarrow \rr (\ep_1 m, \ep_2 n) \in \R^3$$
such that $\Delta_1 \rr$ and $\Delta_2 \rr$ are linearly independent for any $m,n$, where  $\Delta_j$ is defined by
\be  \label{Delta}
   \Delta_j f =  \frac{T_j f - f}{\ep_j} \ ,
\ee
and $f : \ep_1 \Z \times \ep_2 Z \rightarrow \R^3$. 
 In other words, we consider the case $\T_1 = \ep_1 \Z$, $\T_2 = \ep_2 \Z$, where $\ep_1$, $\ep_2$ are fixed constants (the mesh size).  
Therefore, in the discrete case $D_j = \Delta_j$. In particular, for $\ep_1 = \ep_2 =1$ we have $\Delta_j = T_j - 1$.

The discrete analogue of pseudospherical surfaces endowed with asymptotic Chebyshev coordinates is defined as follows (compare \cite{Sau2,Wun}). Weak Chebyshev coordinates were discretized in a similar way.

\begin{Def}[\cite{BP-pseudo}] \label{pseudodis} 
 Discrete asymptotic weak Chebyshev net (discrete $K$-surface)  
is an  immersion $\rr: \ep_1 \Z \times \ep_2 \Z \rightarrow \R^3$ such that for any $m,n$
\begin{itemize}
\item $\Delta_1 \rr \cdot \Delta_1 \rr = E (m)$ , \ $\Delta_2 \rr \cdot \Delta_2 \rr = G (n)$  , \ (weak Chebyshev net)  , discrete Chebyshev nets correspond to $E = G = 1$, 
\item the points $\rr$, $T_1 \rr$, $T_2 \rr$, $T_1^{-1} \rr$, $T_2^{-1} \rr$  are  coplanar  \ (asymptotic net).
\end{itemize}
\end{Def}

The plane containing $\rr$, $T_1 \rr$, $T_2 \rr$, $T_1^{-1} \rr$, $T_2^{-1} \rr$  can be interpreted as the discrete analogue of the tangent plane and   
\be
  \nn := \frac{ \Delta_1 \rr \times \Delta_2 \rr }{ | \Delta_1 \rr \times \Delta_2 \rr | } = \frac{ \Delta_1 \rr \times \Delta_2 \rr }{\sqrt{ (\Delta_1 \rr)^2 (\Delta_2 \rr)^2  - (\Delta_1 \rr \cdot \Delta_2 \rr)^2 } } \ ,
\ee
is the discrete analogue of the normal vector (here the cross means the vector product).

\section{Some old results in a new form}

In order to obtain the explicit similarity between smooth and discrete cases we will reformulate the definition of discrete asymptotic nets and derive another formula for the discrete Gaussian curvature of asymptotic Chebyshev nets. 

\begin{prop}  \label{old1}
For any discrete immersion $\rr$
\[
\Delta_1 \nn \cdot \Delta_1 \rr = 0 \quad \Longleftrightarrow \quad \Delta_1 \rr \ , \ T_1 (\Delta_1 \rr)\ , \ T_1 (\Delta_2 \rr) \ {\rm are \ coplanar}  .
\]
\[
\Delta_2 \nn \cdot \Delta_2 \rr = 0 \quad \Longleftrightarrow \quad \Delta_2 \rr \ , \ T_2 (\Delta_1 \rr)\ , \ T_2 (\Delta_2 \rr) \ {\rm are \ coplanar}  .
\]
\end{prop}

\begin{Proof} From the definition of $\nn$ it follows: $\nn \cdot \Delta_1 \rr = 0$, $T_1 \nn \cdot T_1 \Delta_1 \rr = 0$ and $T_1 \nn \cdot T_1 \Delta_2 \rr = 0$. Then $\Delta_1 \nn \cdot \Delta_1 \rr = 0  \ \Longleftrightarrow \ T_1 \nn \cdot \Delta_1 \rr = \nn \cdot \Delta_1 \rr $. Hence, $T_1 \nn \cdot \Delta_1 \rr = 0$. Therefore, 
$\Delta_1 \rr$, $T_1 \Delta_1 \rr$ and $T_1 \Delta_2 \rr$ are co-planar.
The proof of the second statement is similar.
\end{Proof}

\begin{cor} \label{cordis} For any discrete immersion the points 
$\rr$,   $T_1 \rr$, $T_2 \rr$, $T_1^{-1} \rr$, $T_2^{-1} \rr$ \ are  coplanar if and only if \ 
$\Delta_1 \nn \cdot \Delta_1 \rr = 0$ and  
$\Delta_2 \nn \cdot \Delta_2 \rr = 0$. In other words, a discrete immersion $\rr : \ep_1 \Z \times \ep_2 \Z \rightarrow \R^3$ is asymptotic iff
\be
\Delta_1 \rr \cdot \Delta_1 \nn = \Delta_2 \rr \cdot \Delta_2 \nn = 0 \ , 
\ee 
which is a discrete analogue of \rf{asym}. 
\end{cor}

\begin{prop} \label{Kprop}
For any discrete asymptotic weak Chebyshev net, $K$ defined by 
\be  \label{Kdisc}
 K : = - \frac{ (\Delta_1 \nn\cdot \Delta_2 \rr ) (\Delta_2 \nn \cdot \Delta_1 \rr )   }{ (\Delta_1 \rr)^2 (\Delta_2 \rr)^2 - (\Delta_1 \rr \cdot \Delta_2 \rr )^2 }  
\ee
is constant (i.e., does not depend on $m,n$).  
\end{prop}

\begin{Proof} We consider the tetrahedron $ABCD$: 
$\rr \equiv A$, $T_1 \rr \equiv B$, $T_2 \rr \equiv D$, $T_1 T_2 \rr \equiv C$. Taking into account Definition~\ref{pseudodis}, we have 
\be  \label{leng}
|\vec{AB}| = |\vec{DC}| = \ep_1 | \Delta_1 \rr | \ , \quad  |\vec{AD}| = |\vec{BC}| = \ep_2 | \Delta_2 \rr | . 
\ee
We denote by $h_{AB}^D$ the height of the triangle $ABD$ perpendicular to $AB$, and by $H^D$ the height of the tetrahedron $ABCD$ perpendicular to the base $ABC$, etc. Then $\theta_1$ denotes the angle between $\nn$ and $T_1 \nn$ (i.e., beteween the planes $ABC$ and $ABD$) and $\theta_2$ denotes the angle between $\nn$ and $T_2 \nn$ (i.e., between $ABD$ and $ACD$). Note that the angle between $ABC$ and $BCD$ is $T_1 \theta_2$, and the angle between $ACD$ and $BCD$ is $T_2 \theta_1$. Finally, $\phi$ is the angle between $\Delta_1 \rr$ and $\Delta_2 \rr$, i.e.,
\be
\Delta_1 \rr \cdot \Delta_2 \rr = |\Delta_1 \rr| |\Delta_2 \rr| \cos\phi \ . 
\ee
From elementary geometric considerations we have:
\be \ba{l}  \label{element}
H^B = h_{AD}^B \sin\theta_2 \ , \quad  h_{AD}^B = |\vec{AB}| \sin\phi \ , \quad h_{CD}^B = |\vec{BC}| \sin\phi \ ,
\\[2ex]
H^D = h_{AB}^D \sin\theta_1 \ , \quad 
h_{AB}^D = |\vec{AD}| \sin\phi \ , \quad h_{BC}^D = |\vec{DC}| \sin\phi \ , \\[2ex]
H^B = h_{CD}^B \sin T_2\theta_1 \ , \quad H^D = h_{BC}^D \sin T_1\theta_2 \ , \quad H^D = H^B \ .
\ea \ee
The last equation results from the comparison of two formulae for the volume of the tetrahedron: 
$H^D P_{ABC}  = H^B P_{ACD}$, where $P_{ABC} = P_{ACD}$ because the triangles $ABC$ and $ACD$ are congruent. From  \rf{element} we obtain: 
\[
\frac{\sin\theta_1}{| \vec{AB} |} = \frac{\sin\theta_2}{ | \vec{AD} |} = \frac{\sin T_2\theta_1}{ | \vec{DC} | } = \frac{\sin T_1 \theta_2}{ | \vec{BC} | } ,
\]
which implies
\be  \label{tors}
\frac{\sin\theta_1}{\ep_1 |\Delta_1 \rr | } = \frac{\sin\theta_2}{\ep_2 | \Delta_2 \rr | } = \const \ .
\ee
Then,
\be \ba{l} \dis   \label{main}
\Delta_1 \nn \cdot \Delta_2 \rr = \frac{T_1 \nn \cdot \Delta_2 \rr}{\ep_1} = \frac{| \Delta_2 \rr | }{\ep_1} \frac{H^D}{|\vec{AD}|} = \frac{| \Delta_2 \rr | \sin\theta_1 \, \sin\phi}{\ep_1} \ , \\[3ex] \dis
\Delta_2 \nn \cdot \Delta_1 \rr = \frac{T_2 \nn \cdot \Delta_1 \rr}{\ep_2} = \frac{| \Delta_1 \rr | }{\ep_2} \frac{H^B}{|\vec{AB}|} = \frac{| \Delta_1 \rr | \sin\theta_2 \, \sin\phi }{\ep_2} \ , \\[4ex]
(\Delta_1 \rr)^2 (\Delta_2 \rr)^2 - (\Delta_1 \rr \cdot \Delta_2 \rr )^2 = 
(\Delta_1 \rr)^2 (\Delta_2 \rr)^2 \sin^2 \phi \ .
\ea\ee
Therefore, computing \rf{Kdisc}, we obtain
\be
K = - \frac{\sin\theta_1 \sin\theta_2}{\ep_1 \ep_2 |\Delta_1 \rr | |\Delta_2 \rr | }  \ ,
\ee
and, taking into account \rf{tors}, we complete the proof.
\end{Proof}

$K$ given by the formula \rf{Kdisc} can be considered as a natural discrete analogue of the Gaussian curvature \rf{Kcontas}. 
Wunderlich \cite{Wun}, in the case of discrete Chebyshev nets ($\theta_1 = \theta_2 = \theta$, $|\Delta_1 \rr| = |\Delta_2 \rr| = 1$ and $\ep_1 = \ep_2 = \ep$), proposed a similar definition: 
\be
K' = - \frac{\sin^2\theta}{\ep^2 \cos\theta} \ .
\ee
Because in this case $\theta =\const$ (compare \rf{tors}), then, obviously, both $K$ and $K'$ are constant. 
In the continuous limit $\theta \rightarrow 0$ which implies $K' \rightarrow K$.

\section{Pseudospherical surfaces on time scales}

Corollary~\ref{cordis} shows that the assumptions of Definition~\ref{pseudodis} can be expressed completely in terms of  delta derivatives. 
Therefore, the extension of this definition on arbitrary time scales is straightforward. First, given an immersion $\rr$ on a time scale, we define the normal vector
\be  \label{timenor}
\nn := \frac{ D_1 \rr \times D_2 \rr }{| D_1 \rr \times D_2 \rr | }  \ .
\ee

\begin{Def} \label{pseudotime} An immersion $\rr: \T_1 \times \T_2 \ni (t_1, t_2) \rightarrow \rr (t_1, t_2) \in \R^3$ such that for any $t \equiv (t_1, t_2) \in \T_1 \times \T_2$
\begin{itemize}
\item $\rr$ is completely delta differentiable ,
\item $\nn$ is completely delta differentiable ,
\item $(D_1 \rr)^2 = E (t_1)$, \ $(D_2 \rr)^2  = G (t_2)$ \ ,
\item $D_1 \nn \cdot D_1 \rr =  D_2 \nn \cdot D_2 \rr = 0$ \ ,
\end{itemize}
is called an asymptotic weak Chebyshev net on the time scale $\T \equiv \T_1 \times \T_2$ (or, in particular case $E = G =1$, an asymptotic Chebyshev net). 
\end{Def}

In the continuous and discrete cases asymptotic weak Chebyshev nets have constant negative Gaussian curvature (see Propositions~\ref{czeb} and \ref{Kprop}) and, as a consequence, they can be identified with pseudospherical surfaces. This is true also in the general case.  

\begin{Th}  \label{KTh}
For any asymptotic Chebyshev net on a time scale $\T = \T_1 \times \T_2$,  $K$ defined by
\be  \label{Ktime}
 K = - \frac{ (D_1 \nn\cdot D_2 \rr ) (D_2 \nn \cdot D_1 \rr )   }{ (D_1 \rr)^2 (D_2 \rr)^2 - (D_1 \rr \cdot D_2 \rr )^2 }  
\ee
is constant. 
\end{Th}

\begin{Proof}
It is sufficient to show that $D_1 K = D_2 K = 0$ at any $t \in \T$. 
If $t$ is both 1-right-dense and 2-right-dense, we repeat the standard proof of Proposition~\ref{czeb}. Namely, using Codazzi equations (resulting from compatibility conditions, i.e., $\nn \cdot \rr,_{jjk} = \nn \cdot \rr,_{jkj}$) we show that 
\[
k,_1 =  k,_2 = 0 \ , \quad  k = \frac{M}{\sqrt{E} \sqrt{G} \sin\phi} \ .
\]
The formula \rf{Ktime}  yields $K = - k^2$,  compare \rf{Kcontas}. Hence $K,_1  =  K,_2 = 0$. 

If $t$ is both 1-right-scattered and 2-right-scattered, we use the proof of Proposition~\ref{Kprop}. We point out, however, that the proof of Proposition~\ref{old1} (which is crucial in order to identify sides of the tetrahedron with appropriate tangent planes) needs a modification. We have to use the assumption about complete delta differentiability of $\rr$. Indeed, if (for instance) $T_1 t$ is 1-right-dense, then without this assumption $T_1 (\Delta_2 \rr)$ does not have to be perpendicular to $T_1 \nn$.

In the ``mixed'' case the proof is also straightforward (although it seems to be most cumbersome). Let, for instance, $t$ is 1-right-dense and 2-right-scattered. We use the Frenet basis $\vec\tau$, $\vec\nu$, $\vec\beta$:
\be
\vec\tau = \frac{\rr,_1}{\sqrt{E}} \ , \quad \vec\nu = \frac{\vec\tau,_1}{\kappa \sqrt{E}} \ , \quad \vec\beta = \vec\tau \times \vec\nu \ , 
\ee
where $\kappa$ is the curvature of the line $t_2 = \const$ at $t$. The Serret-Frenet equations read:
\be  \label{Freneteqs}
\tau,_1 = \sqrt{E} \kappa \vec\nu \ , \quad \vec\nu,_1 = \sqrt{E} ( {\tilde \kappa} \vec\beta - \kappa \vec\tau ) \ , \quad \vec\beta,_1 = - \sqrt{E} {\tilde \kappa} \vec\nu \ ,
\ee
where $\tilde \kappa$ is the second curvature (or the torsion). 
We define a unit vector $\vec d$
\be
  \vec d := \frac{D_2 \rr}{\sqrt{G}} \ , \quad \nn = \frac{\vec\tau \times \vec d}{\sin\phi} \ .
\ee
From $\nn,_1 \cdot \rr,_1 = 0$ we derive $(\rr,_1 \times \rr,_{11}) \cdot {\vec d} = 0$. Hence $\vec d \perp \vec\beta$. Then 
\be  \label{dvec}
\vec d = \vec\tau  \cos\phi + \vec\nu \sin\phi \ , \quad \nn = \beta \ .
\ee
$D_2 \rr \cdot D_2 \rr = G (t_2)$ implies $\vec d \cdot T_2 \rr,_1 = \rr,_1 \cdot \vec d$, and $D_2 \nn \cdot D_2 \rr = 0$ implies $T_2 \nn \cdot \vec d = 0$. Then, from $T_2 \nn,_1 \cdot T_2 \rr,_1 = 0$ we get $T_2 \nn \cdot T_2 \rr,_{11} = 0$, i.e., $T_2 \nn = T_2 \vec\beta$. Hence
\be
\vec d = T_2 \vec\tau \cos\phi + T_2 \vec\nu \sin\phi \ .
\ee
Therefore, introducing an additional angle $\vartheta$ and performing two rotations, we can express the basis $T_2 \vec\tau$, $T_2 \vec\nu$, $T_2 \vec\beta$ as follows: 
\be  \label{T2} \ba{l}
T_2 \vec\tau = \vec\tau (\cos^2\phi + \cos\vartheta \sin^2\phi) + \vec\nu \sin\phi \cos\phi (1 - \cos\vartheta) + \vec\beta \sin\phi \sin\vartheta \ ,
\\[2ex]
T_2 \vec\nu = \vec\tau \sin\phi \cos\phi (1 - \cos\vartheta) + \vec\nu (\sin^2\phi + \cos\vartheta \cos^2\phi) - \vec\beta \cos\phi \sin\vartheta \ ,
\\[2ex]
T_2 \vec\beta = - \vec\tau \sin\phi \sin\vartheta + \vec\nu \cos\phi\sin\vartheta + \vec\beta \cos\vartheta \ .
\ea \ee 
On the other hand, we have $T_2 \rr = \rr + \ep {\vec d} \sqrt{G}$, where $\ep = \sigma (t_2) - t_2$. Differentiating it (remember that $G,_1 = 0$) and using \rf{dvec}, \rf{Freneteqs} we get
\be  \label{T2r1}
\frac{T_2 \vec\tau - \vec\tau}{ \ep \sqrt{G} } = \left( \kappa + \frac{\phi,_1}{\sqrt{E}} \right) \left( - \vec\tau \sin\phi + \vec\nu \cos\phi \right) + {\tilde \kappa} \vec\beta \sin\phi \ . 
\ee
Comparing \rf{T2} with \rf{T2r1} we explicitly express $\kappa$ and $\tilde\kappa$ by $\phi$ and $\vartheta$: 
\be  \label{kappy}
\ep \tilde\kappa \sqrt{G} = \sin\vartheta \ , \quad \ep \sqrt{G} \left( \kappa + \frac{\phi,_1 }{\sqrt{E} } \right) = (1 - \cos\vartheta) \sin\phi \ . 
\ee
Substituting \rf{kappy} to the compatibility conditions of equations \rf{Freneteqs} and \rf{T2} we get 
\be  \label{stale}
\vartheta,_1 = 0 \ , \quad  T_2 \left( \frac{\sin\vartheta}{\ep \sqrt{G}} \right) = \frac{\sin\vartheta}{\ep \sqrt{G}} \ ,
\ee
Taking into account $\nn = \vec\beta$ and equations \rf{Freneteqs}, \rf{T2}, \rf{kappy} we compute $K$ using the formula \rf{Ktime}: 
\be
K = - \frac{ (T_2 \vec\beta \cdot \rr,_1)(D_2 \rr \cdot \beta,_1)}{\ep E G \sin^2 \phi} = - \frac{\tilde\kappa \sin\vartheta}{\ep \sqrt{G}} = - \frac{\sin^2\vartheta}{\ep^2 G} \ . 
\ee
Hence, by virtue of \rf{stale} and because $\ep, G$ by assumption do not depend on $t_1$, we have $K,_1 = 0$ and $D_2 K = 0$. 
\end{Proof}

\section{The Lax pair and the Sym formula}   \label{Lax}

Since a pioneering work of Sym (\cite{S-2}, see also \cite{Sym})  smooth pseudospherical surfaces can be constructed from solutions of the corresponding spectral problem (Lax pair) using the so called Sym formula $\Psi^{-1} \Psi,_\lambda$. This approach was extended on discrete surfaces by Bobenko and Pinkall \cite{BP-pseudo,BP-rev}. 

The results of \cite{Ci-hyper} show that relatively weak assumptions on the spectral problem yield smooth pseudospherical surfaces in asymptotic coordinates (asymptotic weak Chebyshev nets). Motivated by these results   
we consider the following system of quaternion-valued linear partial differential equations (the Lax pair) on a time scale ${\mathbb T}_1 \times {\mathbb T}_2$
\be  \ba{l}   \label{problin}
D_1 \Psi = U \Psi \ , \quad U = \lambda (a \e_1 + b \e_2) + c \e_3 + h \ , \\[2ex]
D_2 \Psi = V \Psi \ , \quad V = \lambda^{-1} (p \e_1 + q \e_2) + r \e_3 + s
\ea \ee
where $a, b, c, h, p, q, r, s$ are real functions on 
${\mathbb T}_1 \times {\mathbb T}_2$. Thus (for real $\lambda$) $U, V$ take values in $\mathbb H$, 
and, as a consequence, $\Psi$ is also ${\mathbb H}$-valued. 

The compatibility conditions yield the following system of nonlinear equations:
\be  \label{cc}
D_2 U - D_1 V + \sigma_2 (U) V - \sigma_1 (V) U  = 0 
\ee
(we recall that here $\sigma_1$, $\sigma_2$ denote jump operators, not to be confused with Pauli matrices).

Given $\Psi$ satisfying a Lax pair of the form $D_1 \Psi = U \Psi$, $D_2 \Psi = V \Psi$,  we define an immersion 
${\mathbf r}: {\mathbb T}_1 \times {\mathbb T}_2 \rightarrow {\rm Im\, } {\mathbb H} \simeq {\mathbb E}^3$ by the (modified) Sym formula
\be  \label{Sym}
{\mathbf  r} = \Pi ( \Psi^{-1} \Psi,_\lambda ) \ ,
\ee
where $\Pi$ is the projection \rf{Pi}. Using \rf{Leib} we compute
\[
D_j {\mathbf  r} = \Pi ( - \sigma_j (\Psi^{-1}) (D_j \Psi) \Psi^{-1} \Psi,_\lambda + \sigma_j (\Psi^{-1}) ( U_j,_\lambda \Psi + U_j \Psi,_\lambda ) ) \ ,
\]
where $U_1 := U$, $U_2 := V$. Hence
\be \ba{l} \label{Dr}
D_1 {\mathbf  r} = \Pi ( (\sigma_1 (\Psi))^{-1} U,_\lambda \Psi) \ , \\[2ex]
D_2 {\mathbf  r} = \Pi ( (\sigma_2 (\Psi))^{-1} V,_\lambda \Psi) \ .
\ea \ee

\begin{Th}  \label{Ksym}
Let ${\mathbf  r} : \T_1 \times \T_2 \rightarrow \R^3$  is the surface defined by \rf{Sym}, where $\Psi$ satisfies the Lax pair \rf{problin}. Then the coordinates $t_1, t_2$ are asymptotic, and the formula \rf{Ktime} yields a constant value $K = - 4 \lambda^2$.   
\end{Th}

\begin{Proof}
We will check separately right-dense points and right-scattered points. At $j$-right-dense points $\sigma_j (\Psi) = \Psi$ and
\be  
D_1 {\mathbf  r} = \Psi^{-1} (a \e_1 + b \e_2) \Psi \ , \qquad
D_2 {\mathbf  r} = - \lambda^{-2} \Psi^{-1} (p \e_1 + q \e_2) \Psi \ , 
\ee
while at $j$-right-scattered points $\sigma_j (\Psi) = (1 + \ep_j U_j ) \Psi$ and
\be  \ba{l} \displaystyle
D_1 {\mathbf  r} = \Psi^{-1} \left( \frac{(a + \ep_1 a h + \ep_1 bc) \e_1 + (b + \ep_1 b h - \ep_1 a c ) \e_2 }{(1 + \ep_1 h )^2 + \ep_1^2 c^2 + \ep_1^2 \lambda^2 (a^2 + b^2) } \right) \Psi \ , \\[4ex] \dis
D_2 {\mathbf  r} = -  \Psi^{-1} \left( \frac{(p + \ep_2 p s + \ep_2 q r) \e_1 + (q + \ep_2 q s - \ep_2 p r ) \e_2 }{ \lambda^2 (1 + \ep_2 s )^2 + \lambda^2 \ep_2^2 r^2 + \ep_2^2 (p^2 + q^2) } \right) \Psi \ , 
\ea \ee
In any case the normal vector (compare \rf{timenor}) can be chosen as  
\be
  {\mathbf  n} = \Psi^{-1} \e_3 \Psi \ .
\ee
At $j$-right-dense points $D_j {\mathbf n} = \Psi^{-1} [ \e_3, U_j ] \Psi$. Therefore 
\be
 D_1 {\mathbf n} = 2 \lambda \Psi^{-1} (a \e_2 - b \e_1 ) \Psi \ , \qquad D_2 {\mathbf n} = 2 \lambda^{-1} \Psi^{-1} (p \e_2 - q \e_1 ) \Psi \ . 
\ee
At $j$-right-scattered points $\ep_j D_j {\mathbf n} = (\sigma_j (\Psi))^{-1} \e_3 \sigma_j (\Psi) - \Psi^{-1} \e_3 \Psi$, hence, after straightforward computations
\be  \ba{l} \dis
D_1 {\mathbf n} =  2 \lambda \Psi^{-1} \left(  \frac{ \ep_1 c (a \e_1 + b \e_2) -  (1 + \ep_1 h) ( b \e_1 - a \e_2) + C_1 \e_3 }{ (1 + \ep_1 h )^2 + \ep_1^2 c^2 + \lambda^2 \ep_1^2 (a^2 + b^2) } \right) \Psi \ ,
\\[4ex]  \dis
D_2 {\mathbf n} = 2 \lambda \Psi^{-1} \left(  \frac{ \ep_2 r (p \e_1 + q \e_2) -   (1 + \ep_2 s) ( q \e_1 - p \e_2) + C_2 \e_3 }{\lambda^2 (1 + \ep_2 s )^2 + \lambda^2 \ep_2^2 r^2 + \ep_2^2 (p^2 + q^2) } \right) \Psi \ ,
\ea \ee
where 
\[ \ba{l} \dis 
2 \lambda C_1 = 2 h + \ep_1 h^2 + \ep_1 c^2 - \lambda^2 \ep_1 (a^2 + b^2)  \ , \\[2ex] \dis
2 \lambda C_2 = \lambda^2 ( 2 s + \ep_2 s^2 + \ep_2 r^2) - \ep_2 (p^2 + q^2)  \ . 
\ea \]
We check that $D_1 {\mathbf n} \cdot D_1 {\mathbf r} =  D_2 {\mathbf n} \cdot D_2 {\mathbf r} = 0$ and (after cumbersome computations) 
\be
    \frac{ (D_1 {\mathbf n} \cdot D_2 {\mathbf r} ) (D_2 {\mathbf n} \cdot D_1 {\mathbf r} )   }{ (D_1 {\mathbf r} )^2 (D_2 {\mathbf r} )^2 - (D_1 {\mathbf r} \cdot D_2 {\mathbf r} )^2 }  = 4 \lambda^2 \ ,
\ee
which ends the proof. 
The result is the same for points of any kind (right-dense or right-scattered in one or both directions)! 
\end{Proof}

\section{The Darboux-B\"acklund transformation}

The standard Zakharov-Shabat construction of the Darboux matrix (see, for instance, \cite{Ci-dbt}) can be extended on arbitrary time scales. 
We consider the transformation $\tilde \Psi = B \Psi$ (where $B$ is the Darboux matrix). Then 
\be  \ba{l}
\tilde U = D_1 (B) B^{-1} + \sigma_1 (B) U B^{-1} \ , \\[2ex]
\tilde V = D_2 (B) B^{-1} + \sigma_2 (B) V B^{-1} \ . 
\ea \ee
We confine ourselves to the simplest Darboux matrix $B$ such that 
\be  
  B = N \left( 1 + \frac{\lambda_1 - \mu_1}{\lambda - \lambda_1} P \right) \ , 
\quad B^{-1} = \left( 1 + \frac{\mu_1 - \lambda_1}{\lambda - \mu_1} P \right) N^{-1} \ ,
\ee
where $P^2 = P$. The projector $P$ has to satisfy the system
\be  \ba{l}  \label{system}
D_1 (P) (1 - P) + \sigma_1 (P) U (\lambda_1) (1 - P) = 0 \ , \\[2ex]
D_2 (P) (1 - P) + \sigma_2 (P) V (\lambda_1) (1 - P) = 0 \ , \\[2ex]
(I - \sigma_1 (P)) ( - D_1 P + U (\mu_1) P ) = 0 \ , \\[2ex]
(I - \sigma_2 (P)) ( - D_2 P + V (\mu_1) P ) = 0 \ .
\ea \ee
One can show that $P$ given by  
\be
\ker P = \Psi (\lambda_1) {\vec c}_1 \ , \quad {\rm Im}\,P = \Psi (\mu_1) {\vec c}_2 \ ,
\ee
where ${\vec c}_j$ are constant vectors, satisfies \rf{system}. 
Assuming that 
\be
U = u_0 + \lambda u_1  \ , \quad V = v_0 + \frac{1}{\lambda} v_1 \ ,
\ee
we compute the transformation rules for $u_0, u_1, v_0, v_1$:  
\be  \ba{l}
{\tilde u}_1 = \sigma_1 (N) u_1 N^{-1} \ , \\[2ex]
{\tilde u}_0 = (D_1 N) N^{-1} + \sigma_1 (N) \bigg( u_0 + (\lambda_1 - \mu_1) \big( \sigma_1 (P) u_1 - u_1 P \big)  \bigg) N^{-1} \ , \\[3ex]
{\tilde v}_0 = (D_2 N) N^{-1} + \sigma_2 (N) v_0 N^{-1} \ , \\[3ex] \displaystyle
{\tilde v}_1 = \sigma_2 (N) \left( 1 - \frac{\lambda_1 - \mu_1}{\lambda_1} \sigma_2 (P) \right) v_1 \left( 1 - \frac{\mu_1 - \lambda_1}{\mu_1} P \right) N^{-1} \ .
\ea \ee
The properties of the Lax pair (the reduction group):
\be \ba{l} \label{red1}
U (-\lambda) = \e_3 U (\lambda) \e_3^{-1} \ , \\[2ex]
 V (-\lambda) = \e_3 V  (\lambda) \e_3^{-1} \ ,  
\ea \ee
\be \ba{l} \label{red2}
U^\dagger (\bar \lambda) U (\lambda) = \lambda^2 (a^2 + b^2) + c^2 + h^2  \ , \\[2ex]
V^\dagger (\bar \lambda) V (\lambda) = \lambda^{-2} (p^2 + q^2) + r^2 + s^2   \ ,
\ea \ee
impose constraints on the Darboux matrix $B$ (compare \cite{Ci-dbt}): 
\be
  P^\dagger = P \ , \quad   P = \e_3 (1 - P) \e_3^{-1} \ , \quad \lambda_1 = - \mu_1 = i \kappa_1 \quad (\kappa_1 \in \R) \ .
\ee
In particular, $c_2$ and $c_1$ are orthogonal, and $c_2 = \e_3 c_1$. Therefore 
\be
  P = \frac{1}{2} \left( 1 + i {\mathbf p} \right) \ , \quad 
{\mathbf p} := p_1 \e_1 + p_2 \e_2  \ .
\ee
where ${\mathbf p}^2 = -1$, i.e., $ p_1^2 + p_2^2 = 1$. The longest equations of the system \rf{system} simplify 
\be \ba{l}
{\tilde u}_0 = (D_1 N) N^{-1} + \sigma_1 (N) \bigg( u_0 + \kappa_1 \big( u_1 {\mathbf p} - \sigma_1 ({\mathbf p}) u_1 \big)  \bigg) N^{-1} \ , \\[3ex]  \displaystyle
{\tilde v}_1 = \sigma_2 (N) \sigma_2 ({\mathbf p})  v_1 {\mathbf p}^{-1} N^{-1} \ ,
\ea \ee
and the Darboux matrix and its inverse become
\be
   B = \frac{ N ( \lambda - \kappa_1 {\mathbf p} )}{\lambda - i \kappa_1} \ , \qquad 
B^{-1} = \frac{ (\lambda + \kappa_1 {\mathbf p}) N^{-1} }{\lambda + i \kappa_1 } \ .
\ee
Finally, the transformation on the level of surfaces reads
\be  \label{rDBT} 
\tilde {\mathbf r} = {\mathbf r}  + \frac{\kappa_1}{\lambda^2 + \kappa_1^2} \Psi^{-1} {\mathbf p} \Psi \ .
\ee
Therefore, the B\"acklund transformation has exactly the same form as in the continuous and in the discrete case: the segment joining $\tilde {\mathbf r}$ and $\mathbf r$ is tangent to ${\mathbf r}$ and has a constant length. 
The main difficulty (in the case of time scales different from $\R$ or $\ep \Z$) is to find explicit seed solutions.  

\section{Conclusions}

In this paper the notion of pseudospherical immersions is extended on the so called time scales,  unifying the continuous and discrete cases in a single framework. It can be  especially important in the context of the numerical approximation of continuous integrable models.  
Another important problem raised in this paper is a search of possible sets the integrable systems can be considered on.

The Gaussian curvature of discrete pseudospherical surfaces is defined in a way admitting a straightforward extension on time scales (Proposition~\ref{Kprop}). Surprisingly, the simple formula \rf{Ktime} turns out to be valid for pseudospherical surfaces in asymptotic coordinates on any time scales (Theorem~\ref{KTh}). 
The range of its applicability will be further investigated. 

The quaternion-valued spectral problem \rf{problin} for pseudospherical surfaces in asymptotic coordinates has  very general form. Actually, Theorem~\ref{Ksym} generalizes some results (isospectral case) of my earlier paper \cite{Ci-hyper}  not only on the discrete case, but on arbitrary time scales. The Darboux-B\"acklund transformation \rf{rDBT} can be used to generate explicit pseudospherical surfaces (soliton solutions) on some interesting, non-standard, time scales. The work in this direction is in progress. 

It would be interesting to extend any other results of the integrable discrete geometry on arbitrary time scales.

{\it Acknowledgements.} I am grateful to Iwona \'Swis\l ocka for a cooperation \cite{Sw-mgr}, to Zbigniew Hasiewicz for helpful discussions, and to Klara Janglajew for turning my attention on references \cite{BG,Hi,Hi2}. My work was partially supported by the Polish Ministry of Science and Higher Education (grant No.\  1 P03B 017 28).

\end{document}